# PRODUCTS OF SUBSETS OF GROUP THAT EQUAL THE GROUP
Oleksandr Vyshnevetskiy

## Abstract

Let $P$ be a probability on a finite group $G$, $P^{(n)}$ $n$-fold convolution of $P$ on $G$. Under mild condition, $P^{(n)}$ at $n \to \infty$ converges to the uniform probability on the group $G$. If $A = \{g \in G, P(g) \neq 0\}$ be the carrier of the probability $P$, then $A^n = \{a_1 \cdot \ldots \cdot a_n, \ a_1, \ldots, a_n \in A\}$ be the carrier of probability $P^{(n)}$.

One of necessary and sufficient conditions for the mentioned convergence is: sequence $A^n$ ($n \to \infty$) stabilizes on $G$, i.e. $A^k = A^{k+1} = \ldots = G$ for a natural number $k$. In other words, product of some multipliers equal to $A$ is $G$. The carrier $A$ is in general case any nonempty subset of group $G$. In the paper we find a condition under which product of some subsets of $G$ is $G$.

## Introduction

Let $G$ be a finite group, $P^{(n)}$ $n$-fold convolution of a probability $P$ on the group $G$. Convergence of $P^{(n)}$ to the uniform probability $U(g) = |G|^{-1}$ $(g \in G)$ of the group $G$ is considered in many papers (see e.g. survey [1] or [2] for groups which are countable sums of finite groups).

Let $A = \{g \in G, P(g) \neq 0\}$ be the carrier of $P$. It was proved in [3] that $P^{(n)}$ converges to $U(g)$ if and only if sequence $A^n$ ($n \to \infty$) stabilizes on $G$, i.e. $A^k = A^{k+1} = \ldots = G$ for a natural number $k$. In other words, product of some multipliers equal to $A$ is $G$. The carrier $A$ can be in general any nonempty subset of group $G$. In the paper we find a condition under which product of some subsets of $G$ is $G$. The first condition for two subsets was proved in [4]. One also can consider such product as a factorization of $G$.

For subsets $A, B \subset G$ their product is defined as $A \cdot B = \{a \cdot b, a \in A, b \in B\}$, complement $\overline{A} = G \setminus A$ and $A^{-1} = \{a^{-1}, \ a \in A\}$. For two subsets of $G$ one of the mentioned conditions for their product is

<u>Lemma 1</u> [4]. If $|A| + |B| > |G|$, then $A \cdot B = G$.

The condition is sufficient but not necessary (see example 4 below, where $|A| + |B| < |G|$).

<u>Corollary</u>. If $|A| > \frac{1}{2}|G|$, then $G = A \cdot A$ and $G = A \cdot A^{-1}$.

Note that one cannot replace the inequality with the equality $|A| = \frac{1}{2}|G|$. Indeed,

for a subgroup $A$, $|G:A|=2$, we have $A \cdot A = A \neq G$ and $AA^{-1} = A^{-1}A = A \neq G$.

## Main result

We generalize Lemma 1 to the case of an arbitrary number $n \geq 2$ of nonempty subsets $A_1,...,A_n \subset G$. Let $Z$ be the ring of integer numbers, $ZG$ a group algebra of group $G$ over the ring $Z$. For an arbitrary subset $X$ in $G$ let

$$[X] = \sum_{x \in X} x \in ZG, \quad [\varnothing] = 0$$

Since $g[G]=[G]g=[G]$ for any $g \in G$, then

$$[X][G]=[G][X]=|X|[G] \qquad (1)$$

For nonempty subsets $A_1,...,A_n \subset G$ let $B=(A_1,...,A_n)$, $\overline{B}=(\overline{A}_1,...,\overline{A}_n)$,

$$N_B(g) = \left|\{(a_1,...,a_n) \in A_1 \times ... \times A_n\}, a_1 \cdot ... \cdot a_n = g\right|$$

$$N_{\overline{B}}(g) = \left|\{(\overline{a}_1,...,\overline{a}_n) \in \overline{A}_1 \times ... \times \overline{A}_n\}, \overline{a}_1 \cdot ... \cdot \overline{a}_n = g\right|$$

<u>Theorem 2</u>. For any element $g \in G$ we have

$$N_B(g) - (-1)^n N_{\overline{B}}(g) = d(B) \qquad (2)$$

where number

$$d(B) = \frac{1}{|G|}\left(\prod_i |A_i| - (-1)^n \prod_i |\overline{A}_i|\right) \qquad (3)$$

(here and up to the end of the proof, in all products and sums, summation indices change by default from 1 to $n$).

Proof. As $[A_i] = [G \setminus \overline{A}_i] = [G] - [\overline{A}_i]$ and $[G][\overline{A}_i] = [\overline{A}_i][G]$, then

$$\prod_i [A_i] = ([G]-[\overline{A}_1])...([G]-[\overline{A}_n]) = [G]^n - \sum_i [\overline{A}_i][G]^{n-1} +$$

$$+ \sum_{i<j} [\overline{A}_i][\overline{A}_j][G]^{n-2} - ... + (-1)^{n-1} \sum_{i_1<...<i_{n-1}} [\overline{A}_{i_1}]...[\overline{A}_{i_{n-1}}][G] + (-1)^n \prod_i [\overline{A}_i]$$

From (1) we have $[G]^k = [G]^{k-1}[G] = |G|^{k-1}[G]$ for $k=1,...n$, so the last equality can be written as

$$\prod_i [A_i] - (-1)^n \prod_i [\overline{A}_i] = |G|^{n-1}[G] - |G|^{n-2} \sum_i [\overline{A}_i][G] +$$

$$+ |G|^{n-3} \sum_{i<j} [\overline{A}_i][\overline{A}_j][G] - ... + (-1)^{n-1} \sum_{i_1<...<i_{n-1}} [\overline{A}_{i_1}]...[\overline{A}_{i_{n-1}}][G] =$$

$$= \left(|G|^{n-1} - |G|^{n-2}\sum_i [\overline{A}_i] + |G|^{n-3}\sum_{i<j}[\overline{A}_i][\overline{A}_j] + (-1)^{n-1}\sum_{i_1<...<i_{n-1}}[\overline{A}_{i_1}]...[\overline{A}_{i_{n-1}}]\right)[G] =$$

$$= \frac{1}{|G|}\left(|G|^n - |G|^{n-1}\sum_i [\overline{A}_i] + |G|^{n-2}\sum_{i<j}[\overline{A}_i][\overline{A}_j] + ... + (-1)^{n-1}|G|\sum_{i_1<...<i_{n-1}}[\overline{A}_{i_1}]...[\overline{A}_{i_{n-1}}]\right)[G] =$$

$$= \frac{1}{|G|}\left(\prod_i\left((|G|-|\overline{A_i}|)-(-1)^n\prod_i|\overline{A_i}|\right)\right)[G] = \frac{1}{|G|}\left(\prod_i|A_i|-(-1)^n\prod_i|\overline{A_i}|\right)[G].$$

Taking into account (3), we have

$$\prod_i[A_i]-(-1)^n\prod_i[\overline{A_i}] = d(B)[G] \tag{4}$$

As

$$\prod_i[A_i] = \sum_{g\in G} N_B(g)g, \quad \prod_i[\overline{A_i}] = \sum_{g\in G} N_{\overline{B}}(g)g,$$

then equating coefficients in both sides of (4) for each $g\in G$ we get (2).

<u>Corollary 1</u>. $d(B)\in Z$

We can prove the corollary independently. Indeed, since $|\overline{A_i}| = |G|-|A_i| \equiv -|A_i|\pmod{|G|}$ $(i=1,...,n)$, then $(-1)^n\prod_i|\overline{A_i}| \equiv \prod_i|A_i|\pmod{|G|}$, so $|G|$ divides $\prod_i|A_i|-(-1)^n\prod_i|\overline{A_i}|$.

We note that $d(B)$ depends neither on order of subsets $A_1,...,A_n$ in $B$ nor on $g\in G$. So the same is true for $N_B(g)-(-1)^n N_{\overline{B}}(g)$.

<u>Corollary 2.</u> 1) If $n$ is even and $d(B)>0$, then $A_1...A_n = G$.

2) If $n$ is even and $d(B)<0$, then $\overline{A_1}...\overline{A_n} = G$.

3) If $d(B)=0$, then $n$ is even and $A_1...A_n = \overline{A_1}...\overline{A_n}$.

Proof. 1) In this case from (2) $N_B(g) = d(B)+N_{\overline{B}}(g) \geq d(B) > 0$ for any $g\in G$, so $N_B(g)>0$, i.e. $A_1...A_n = G$.

2) If $d(B)<0$, then $d(\overline{B}) = -d(B)>0$, so by point 1), $\overline{A_1}...\overline{A_n} = G$

3) If $d(B)=0$, then from (2) we have $N_B(g) = (-1)^n N_{\overline{B}}(g)$ for any $g\in G$. As $N_B(g)>0$ for some $g\in G$ and $N_{\overline{B}}(g)\geq 0$, then $n$ is even and $N_B(g) = N_{\overline{B}}(g)$ for any $g\in G$. Hence $A_1...A_n = \overline{A_1}...\overline{A_n}$.

In case $n=2$, i.e. for two subsets $A_1, A_2 \subset G$, we get

<u>Corollary 3.</u> 1) If $|A_1|+|A_2|>|G|$, then $A_1 A_2 = G$.

2) If $|A_1|+|A_2|<|G|$, then $\overline{A_1}\,\overline{A_2} = G$.

3) If $|A_1|+|A_2|=|G|$, then $A_1 A_2 = \overline{A_1}\,\overline{A_2}$.

Proof. In case $n=2$ equality (3) has form

$$d(B) = \frac{1}{|G|}(|A_1||A_2|-|\overline{A_1}||\overline{A_2}|) = \frac{1}{|G|}(|A_1||A_2|-(|G|-|\overline{A_1}|)(|G|-|\overline{A_2}|)) = |A_1|+|A_2|-|G|,$$

so points 1)-3) of previous corollary yield points 1)-3) of this one.

For any $A\subset G$ we denote $A^2 = AA$, $A^{-1} = \{a^{-1}, a\in A\}$.

<u>Corollary 4.</u> Let $A,B\subset G$ and $|A|>\frac{1}{2}|G|$. Then

1) If $|B| \geq \frac{1}{2}|G|$, then $A^2 = AB = BA = G$

2) $AA^{-1} = A^{-1}A = G$.

For item 2) we note that $|A| = |A^{-1}|$.

<u>Corollary 5</u>. $A\overline{A} = \overline{A}A$ for any $A \subset G$.

Proof. Since $|A| + |\overline{A}| = |G|$, it follows from point 3) of corollary 3.

In Corollary 3, item 3) (unlike in items 1) and 2)) we proved neither $A_1 A_2 = G$ nor $A_1 A_2 \neq G$. Both the cases are possible if $|A_1| + |A_2| = |G|$, see examples 1 – 3 below.

We note that if $H$ is a subgroup of group $G$, then
$$H \cdot (G \setminus H) \cap H = \varnothing \qquad (5)$$
Indeed, if $h \in H$, $g \in G \setminus H$, then $hg \notin H$ (otherwise $g \in H$).

<u>Example 1</u>. Let $A_1 = H$ be a subgroup of group $G$, $A_2 = G \setminus H$. Then $|A_1| + |A_2| = |G|$. Due to (5), $A_1 A_2 = H \cdot (G \setminus H) \neq G$.

<u>Example 2</u>. Let $G = \langle a, b; a^2 = b^2 = 1, ab = ba \rangle$ be an elementary abelian 2-group of order 4, $A_1 = \langle a, a^2 = 1 \rangle$, $A_2 = \langle b, b^2 = 1 \rangle$ two its subgroups. Again $|A_1| + |A_2| = |G|$, but now $A_1 A_2 = G$.

In example 2 $A_1 \cap A_2 = \{1\} \neq \varnothing$. In next example 3 $|A_1| + |A_2| = |G|$, $A_1 A_2 = G$, but $A_1 \cap A_2 = \varnothing$.

<u>Example 3</u>. Let $G = \langle a, b, c; a^3 = b^3 = c^3 = 1, ab = ba, ac = ca, bc = cb \rangle$ be an elementary abelian 3-group, $A_1 = \{1, a, b\}$, $A_2 = G \setminus H$. Here $A_1 A_2 = G$, i.e. each element $g \in G$ can be written as $g = a_1 a_2$, where $a_1 \in A_1$, $a_2 \in A_2$. Indeed, for $g \in A_2$ the decomposition is $g = 1 \cdot g$ and for elements of $A_1$ it is $1 = a \cdot a^2$, $a = b \cdot b^2 a$, $b = a \cdot a^2 b$. Of course, $A_1 \cap A_2 = \varnothing$.

If we replace $A_1$ with $\overline{A_1}$ and $A_2$ with $\overline{A_2}$ in the examples 1-3, we obtain corresponding examples for $\overline{A_1}$ and $\overline{A_2}$. For instance, from example 2 we obtain case $|\overline{A_1}| + |\overline{A_2}| = |G|$, $\overline{A_1}\,\overline{A_2} = G$, $\overline{A_1} \cap \overline{A_2} \neq \varnothing$.

<u>Example 4</u>. Let $A_1$ be a subgroup of $G$, $|G : A_1| = 2$ and $|G| > 4$. Then $A_1 \triangleleft G$ and $G = A_1 \cup A_1 t$ for any element $t \in G \setminus A_1$. Let $A_2$ be a subset of $G$ such that $\{1, t\} \in A_2$. Then $A_1 A_2 \supset A_1 \cdot \{1, t\} = A_1 \cup A_1 t = G$. So $A_1 A_2 = G$, though $|A_1| + |A_2| < |G|$. So condition $|A_1| + |A_2| > |G|$ is sufficient but not necessary for $A_1 A_2 = G$ (compare with item 1) of corollary 3).

In this example $\overline{A_1}\,\overline{A_2} = G$, because $|A_1| + |A_2| < |G|$. So both the equalities $A_1 A_2 = G$ and $\overline{A_1}\,\overline{A_2} = G$ holds.

Now we generalize Corollary 3 for case of $n \geq 2$ nonempty subsets $A_1,...,A_n \subset G$. Let $A_1 \cdot ... \cdot A_n = C$. We call $A_i \cdot ... \cdot A_j$ $(1 \leq i < j \leq n)$ a subproduct of $C$ of length $j - i$.

We note that $AG = GA = G$ for a nonempty subset $A \subset G$, so if $C$ contains a subproduct equal to $G$, then $C = G$.

<u>Theorem 3</u>. Let $1 \leq i, j \leq n$.

1) If $|A_i| + |A_j| > |G|$ for some $i \neq j$, then $C = G$.

2) If $|A_i| + |A_j| < |G|$ for some $i \neq j$, then $\overline{A}_1 \cdot ... \cdot \overline{A}_n = G$.

3) Otherwise both the cases $C = G$ and $C \neq G$ are possible for any $n > 1$.

Proof. 1) We can consider $j > i$, i.e. $j = i + k$ for some natural number $k$. Induction on length $k$ of subproduct $A_i \cdot ... \cdot A_{i+k}$ in $C$. If $k = 1$, then $A_i \cdot ... \cdot A_{i+k} = A_i \cdot A_{i+1} = G$ by Corollary 3. Product $C$ has a subproduct $A_i \cdot A_{i+1} = G$, so due to the note above $C = G$.

Let $k > 1$. We have
$$A_i \cdot ... \cdot A_{i+k-1} \cdot A_{i+k} \cdot A_{i+k+1} = \left( A_i \cdot ... \cdot A_{i+k-1} \cdot A_{i+k+1} \cdot A_{i+k} \right) \cdot \left[ A_{i+k+1}, A_{i+k} \right] \quad (6)$$
where $\left[ A_{i+k+1}, A_{i+k} \right]$ is commutant of subsets $A_{i+k+1}$ and $A_{i+k}$. The product in parenthesis in (6) is a subproduct of length $k - 1$. By induction hypothesis the product it is $G$. Again product $C$ has a subproduct equal to $G$, so $C = G$.

2) If $|A_i| + |A_j| < |G|$, then
$$|\overline{A}_i| + |\overline{A}_j| = |G| - |A_i| + |G| - |A_j| = 2|G| - (|A_i| + |A_j|) > |G|,$$
so item 2) follows from item 1).

3) If neither of conditions 1), 2) of the theorem is fulfilled, then
$$|A_i| + |A_j| = |G| \quad (1 \leq i, j \leq n) \quad (7)$$

For $n = 2$ example 1 shows a group $G$ with $A_1 A_2 \neq G$ and example 2 a group $G$ with $A_1 A_2 = G$. If $n > 2$, then (7) yields $|A_i| = \frac{|G|}{2}$ $(1 \leq i \leq n)$. In this case:

- if $G$ has a subgroup $A$, $|G : A| = 2$, then $A \cdot A \cdot ... \cdot A = A \neq G$, i.e. $C \neq G$
- if $A_1 \cdot ... \cdot A_{n-1} = G$, then $C = A_1 \cdot ... \cdot A_n = G$ for any subset $A_n$ such that $1 \in A_n$ and $|A_n| = \frac{|G|}{2}$.

O. L. VYSHNEVETSKIY
Department of High Mathematics, Kharkiv National Automobile and Highway University,
61002, Ukraine, Kharkiv, st. Yaroslava Mudrogo, 25
Email address: alexwish50@gmail.com